\documentclass{article}
\usepackage{graphicx}
\usepackage[all]{xy}
\usepackage{tikz}
\usepackage{amsmath}
\usepackage{amsthm}
\usepackage{amssymb}
\begin{document}

\title{Self-Automaton Semigroups}
\author{Alexander McLeman\\
  School of Mathematics and Statistics,\\
  University of St Andrews,\\
  North Haugh,\\
  St Andrews,\\
  KY16 9SS,\\
  United Kingdom\\
  \texttt{alam@st-andrews.ac.uk}}

\maketitle

\begin{abstract}
After reviewing automaton semigroups, we introduce Cayley Automata and the corresponding Cayley Automaton semigroups. We investigate which semigroups are isomorphic to their Cayley Automaton semigroup and give some results for special classes of semigroups. We answer a question posed by Cain relating to the dual construction.
\end{abstract}

\section{Introduction}

Automaton groups are objects that have received considerable attention, most notably from authors such as Grigorchuk, Nekrashevich and Sushchanskii. These are naturally generalised to automaton semigroups, objects which have been studied recently by Cain in \cite{cain09} and by Silva and Steinberg in \cite{silva05}. Automaton semigroups are constructed by considering the actions of the states of the automaton on sequences over an alphabet and generating a semigroup from these induced transformations. An automaton group is an automaton semigroup generated by the bijective induced transformations and their inverses.

There is a natural way of viewing the Cayley table of a finite semigroup as an automaton, which gives rise to the notion of a \emph{Cayley automaton}. Each element of the semigroup corresponds to a state in the automaton and the transitions between states and the outputs from the automaton are defined by right-multiplication in the semigroup. These automata have their foundations in the works of Krohn and Rhodes \cite{krohn65, krohn68}. The automaton semigroups associated to this particular class of automata are termed \emph{Cayley Automaton Semigroups}. These semigroups have been studied in recent years by several authors, such as Maltcev \cite{maltcev09}, Silva and Steinberg \cite{silva05} (who employ wreath products in their work), Mintz \cite{mintz09} (who acts directly on sequences with states of the automaton - the approach we will take in this paper) and Cain \cite{cain09} (who utilises both approaches). This led to some important results, such as determining that all free semigroups of rank at least 2 arise in this way (see Theorem \ref{th1} below) and identifying when a Cayley Automaton semigroup is finite (see Theorem \ref{th2} below). 

It is often natural with expansion-like constructions, of which the Cayley Automaton semigroup construction is an example, to investigate the objects that are invariant under the construction. These objects will be termed \emph{self-automaton semigroups}. This viewpoint is adopted by Cain in \cite[Section 6.4]{cain09} where he gives as examples of self-automaton semigroups semilattices and $I\times I$ rectangular bands with an identity. This led him to pose the following:
\theoremstyle{plain}
\newtheorem{introprob}{Question}[section]
\begin{introprob}
\label{introprob}
Can self-automaton semigroups be classified? The class of such semigroups might consist of precisely those finite bands in which every $\mathcal{D}$-class is square and every topmost $\mathcal{D}$-class is a singleton.
\end{introprob}

At first sight, the requirement for square $\mathcal{D}$-classes may appear surprising. On further inspection of the problem, this may have arisen as a result of the choice by Cain (and also by Maltcev) to act on sequences from the right with states of the automaton, whereas other authors have chosen to adopt the convention of acting from the left. This latter approach of acting from the left appears more natural to us, and in this paper, we will bring Cain's concept of self-automaton semigroups back into this context. 

It is worth emphasising that the outcomes of choosing to act from the left rather than the right go beyond a mere anti-isomorphism. We obtain a richer and more interesting class of semigroups that are self-automaton and are able to interpret fully Cain's notion of self-automaton within our framework.

The main concepts of the paper - those of automaton semigroups and Cayley Automaton semigroups - are introduced formally in sections \ref{sec:autsgp} and \ref{sec:cayley} respectively. In section \ref{sec:self}, self-automaton semigroups are introduced and we explore the initial links between the semigroups and their left-regular representations. Guided by the open problem stated earlier, we next turn our attention to bands in section \ref{sec:bands}, where we prove that a band is self-automaton if and only if its left-regular representation is faithful (see Theorem \ref{th7} below). At this point, it is tempting to rephrase Cain's open problem in our setting of left actions, which would perhaps go as follows:
\theoremstyle{plain}
\newtheorem{quest}[introprob]{Question}
\begin{quest}
\label{quest}
Does the class of self-automaton semigroups consist precisely of those bands which have a faithful left-regular representation?
\end{quest}

In attempting to answer this, we prove positive results for regular semigroups and monoids but show that in general the question has a negative answer i.e. the class of self-automaton semigroups contains more than just bands. Suitable examples of non-band self-automaton semigroups are constructed in section \ref{sec:nonband}. It should be noted that these are not counterexamples to Cain's original problem due to the lack of symmetry in the $\mathcal{D}$-classes. 

We return to Cain's original notion in section \ref{sec:cain} and interpret it fully in our framework before going on to construct an example to address his problem. Finally, further properties of self-automaton semigroups are discussed in section \ref{sec:other} where we prove that for all bands, the arising Cayley Automaton semigroup is itself self-automaton, even if the original band was not. For now, the full classification of self-automaton semigroups appears to be both desirable and elusive. 

\section{Automaton Semigroups}
\label{sec:autsgp}
In defining an automaton, we will follow the definitions given in \cite{cain09} and \cite{silva05}.
\theoremstyle{definition}
\newtheorem{def1}{Definition}[section]
\begin{def1}
\label{def1}
An \emph{automaton} $\mathcal{A}$ is a triple $(Q,B,\delta)$ consisting of a finite set of \emph{states} $Q$, a finite \emph{alphabet} $B$ and a \emph{transition function} $\delta:Q\times B\rightarrow Q\times B.$ We think of the automaton $\mathcal{A}$ as a directed labelled graph with vertex set $Q$ and an edge from $q$ to $r$ labelled by $x|y$ exactly when $\delta(q,x)=(r,y).$ Pictorially we have:
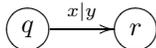
\[\xymatrix{
& *++[o][F]{q}\ar[r]^{x|y}
& *++[o][F]{r}
}\]
\end{def1}
For the purpose of clarity, we make the following definitions which will be adhered to throughout:
\newtheorem{def2}[def1]{Definition}
\begin{def2}
\label{def2}
A \emph{word} is an element of $Q^+$. A \emph{sequence} is an element of $B^*$, consisting of \emph{symbols} from $B.$
\end{def2}
The sequences in $B^*$ are acted on by the states in $Q.$ Throughout this paper, states are defined to act on sequences from the \emph{left}, in contrast to \cite{cain09} and \cite{maltcev09}. We define the action as follows: $q\cdot\alpha$ (the result of state $q$ acting on a sequence $\alpha$) is by definition the sequence outputted by the automaton after starting in state $q$ and reading $\alpha.$ More explicitly, if $\alpha=\alpha_1\alpha_2\ldots\alpha_n$ (for $\alpha_i\in B$) then $q\cdot\alpha=\beta_1\beta_2\ldots\beta_n$ where $\delta(q_{i-1},\alpha_i)=(q_i,\beta_i)$ for $i=1,\ldots,n$ and $q_0=q.$

Notice here that for each symbol read by the automaton exactly one symbol is outputted and so $|q\cdot\alpha|=|\alpha|$ (where $|\alpha|$ denotes the length of the sequence $\alpha$). Such automata are referred to as \emph{synchronous} or \emph{Mealy}  automata, and as such, the action on finite sequences is determined by the action on infinite sequences and vice-versa. An infinite sequence which consists of countably many repetitions of a sequence $\alpha\in B^*$ is denoted by $\alpha^{\omega}.$

We can identify the set of sequences $B^*$ with the regular rooted $|B|-$ary tree, where the vertices are labelled by the elements of $B^*$ (with the root being the empty word). A state $q$ acting on $B^*$ can thus be considered as a transformation of the tree corresponding to $B^*$ which maps a vertex $w$ to the vertex $q\cdot w.$ This action on the tree is adjacency- and level-preserving and is hence an endomorphism of the tree. 

It is natural to extend the action of states to an action of words. A word $q_nq_{n-1}\ldots q_2q_1$ acts on a sequence $\alpha$ as follows:
$$q_n\cdot(q_{n-1}\cdot\ldots(q_2\cdot(q_1\cdot\alpha))\ldots).$$
This gives a homomorphism $\Phi:Q^+\rightarrow \text{End}(B^*)$ where $\text{End}(B^*)$ is the endomorphism semigroup of $B^*.$ We denote the image of $\Phi$ by $\Sigma(\mathcal{A}).$
\theoremstyle{definition}
\newtheorem{def3}[def1]{Definition}
\begin{def3}
\label{def3}
A semigroup $S$ is said to be an \emph{automaton semigroup} if there exists an automaton $\mathcal{A}$ such that $S\cong\Sigma(\mathcal{A}).$
\end{def3}

It would be instructive at this point to give an example of an automaton semigroup being constructed from an automaton. We will follow an example from \cite{cain09}.
\theoremstyle{definition}
\newtheorem{ex1}[def1]{Example}
\begin{ex1}
\label{ex1}
Let $\mathcal{A}=(\{a,b\},\{0,1\},\delta)$ be the automaton below:
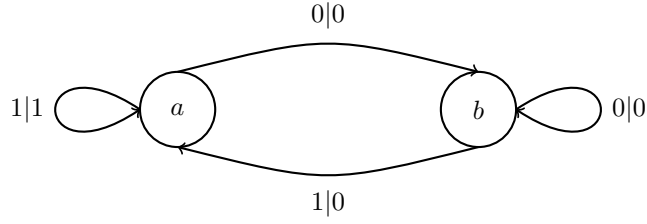
\begin{center}
\begin{tikzpicture}
\draw[thick] (-2,0) circle (0.5cm);
\draw[thick] (2,0) circle (0.5cm); 
\draw[->][thick] (-2.5,0)..controls (-4,1) and (-4,-1)..(-2.5,0);
\draw[->][thick] (2.5,0)..controls (4,1) and (4,-1)..(2.5,0);
\draw[->][thick] (-2,0.5)..controls (0,1)..(2,0.5);
\draw[->][thick] (2,-0.5)..controls (0,-1)..(-2,-0.5);
\draw (-2,0) node {$a$};
\draw (2,0) node {$b$};
\draw (0,1.25) node {$0|0$};
\draw (0,-1.25) node {$1|0$};
\draw (-4,0) node {$1|1$};
\draw (4,0) node {$0|0$};
\end{tikzpicture}
\end{center} 
The transition function $\delta$ is formally defined by:
\begin{center}
\begin{align*}
(a,0)&\mapsto(b,0)\\
(b,0)&\mapsto(b,0)\\
(a,1)&\mapsto(a,1)\\
(b,1)&\mapsto(a,0).\\
\end{align*}
\end{center}

If, for example, $\mathcal{A}$ is in state $a$ and reads the sequence $0011$ then the calculation will proceed as follows:
$a\cdot0011=0(b\cdot011)=00(b\cdot11)=000(a\cdot1)=0001.$

We must consider the actions of $a$ and $b$ on sequences to determine the automaton semigroup. Let $\alpha$ be an infinite sequence. Observe that $b\cdot\alpha$ must start with $0$ for all sequences $\alpha$ and so we may write $b\cdot\alpha=0\beta$ for some sequence $\beta.$ Observe now that $a\cdot b\cdot\alpha=a\cdot 0\beta=0(b\cdot\beta)$ and $b\cdot b\cdot\alpha=b\cdot 0\beta=0(b\cdot\beta).$ Hence in $\Sigma(\mathcal{A})$ we have the relation $ab=b^2.$

We show now that every product in $\Sigma(\mathcal{A})$ can be uniquely expressed as $b^ia^j$ for some $i,j\in\mathbb{N}\cup\{0\}.$ By writing $b^i$ as a state we mean $\underbrace{b\cdot\ldots\cdot b}_{i\text{ times}}.$ Let $i,j\geq 0.$ We have that
\begin{equation}
\label{eqn1}
b^i\cdot a^j\cdot 01^\omega=b^i\cdot 0^{j+1}1^\omega=0^{i+j+1}1^\omega
\end{equation}
and, for $n>i$
\begin{equation}
\label{eqn2}
b^i\cdot a^j\cdot 1^n0^\omega=b^i\cdot 1^n0^\omega=0^i1^{n-i}0^\omega. 
\end{equation}
Now suppose that $b^i\cdot a^j=b^k\cdot a^l$ for some $k,l\geq 0.$ Then by (\ref{eqn2}) $n-i=n-k$ and hence $i=k.$ By (\ref{eqn1}) $i+j+1=k+l+1$ which gives $j=l.$ The semigroup $\Sigma(\mathcal{A})$ is therefore presented by $\langle a,b |  ab=b^2\rangle.$
\end{ex1}

\section{Cayley Automaton Semigroups}
\label{sec:cayley}
In this section we move on to define and discuss the automata and the semigroups arising from them that will form the foundations of the remainder of the paper. The automata are constructed from the Cayley Table of a finite semigroup and are termed \emph{Cayley Automata}. The automaton semigroups arising from these automata are called \emph{Cayley Automaton Semigroups}. They have their origins in \cite{krohn65, krohn68} but more recently have been studied in \cite{cain09, silva05, maltcev09, mintz09}.
\theoremstyle{definition}
\newtheorem{def4}{Definition}[section]
\begin{def4}
\label{def4}
For a finite semigroup $S,$ the \emph{Cayley Automaton} is the automaton $\mathcal{C}(S)=(S,S,\delta)$ where the transition function $\delta$ is defined by $\delta(s,t)=(st,st).$ Note that this is still a synchronous automaton as $st$ is a product in the semigroup $S.$ A typical edge in $\mathcal{C}(S)$ is as follows: 
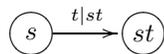
\[\xymatrix{
& *++[o][F]{s}\ar[r]^{t|st}
& *++[o][F]{st}
}\]
Note also that the automaton can be constructed from the right Cayley graph of the semigroup.

A semigroup $T$ is called a \emph{Cayley Automaton semigroup} if there exists a semigroup $S$ such that $T\cong\Sigma(\mathcal{C}(S)).$
\end{def4}
Notice in the definition above that the state set and the alphabet are the same. To avoid any confusion, we will adhere to the following convention:

\newtheorem{not1}[def4]{Notation}
\begin{not1}
\label{not1}
For $s\in S$, the corresponding state in the Cayley Automaton will be written as $\overline{s}$ and the symbol as $s.$
\end{not1}

Let $s\in S$ and $\alpha=\alpha_1\alpha_2\ldots\alpha_n\in S^*.$ Then $$\overline{s}\cdot\alpha=(s\alpha_1)(s\alpha_1\alpha_2)\ldots(s\alpha_1\alpha_2\ldots\alpha_n).$$ With this in mind, we may view each state $\overline{s}$ as a transformation 
$$\overline{s}:\alpha_1\alpha_2\ldots\alpha_n\mapsto(s\alpha_1)(s\alpha_1\alpha_2)\ldots(s\alpha_1\alpha_2\ldots\alpha_n).$$ The semigroup $\Sigma(\mathcal{C}(S))$ is the subsemigroup of $\text{End}(S^*)$ generated by $\{\overline{s}\}_{s\in S}$ under composition of transformations.

Many of the natural questions that arise in the area of Cayley Automaton semigroups concern linking the properties of a semigroup $S$ with the properties of $\Sigma(\mathcal{C}(S)).$ We mention a few of these here before moving on to the question that this paper will address. The following is proved in \cite{silva05}:
\theoremstyle{plain}
\newtheorem{th1}[def4]{Theorem}
\begin{th1}
\label{th1}
Let $G$ be a finite non-trivial group. Then $\Sigma(\mathcal{C}(G))$ is a free semigroup of rank $|G|.$
\end{th1}
Mintz proves the following in \cite{mintz09}, which also appears in \cite{cain09} and \cite{maltcev09}:
\newtheorem{th2}[def4]{Theorem}
\begin{th2}
\label{th2}
Let $S$ be a finite semigroup. Then $\Sigma(\mathcal{C}(S))$ is finite if and only if $S$ is aperiodic.
\end{th2}

\section{Self-Automaton Semigroups}
\label{sec:self}
We turn now to the main theme of the paper and address the question of when a semigroup is isomorphic to its Cayley Automaton semigroup. We make the following definition:
\theoremstyle{definition}
\newtheorem{def5}{Definition}[section]
\begin{def5}
\label{def5}
Let $S$ be a finite semigroup. Then $S$ is \emph{self-automaton} if the map $S\rightarrow\Sigma(\mathcal{C}(S))$ which maps  $s\mapsto\overline{s}$ is an isomorphism.
\end{def5}
Notice that the map $s\mapsto\overline{s}$ is always a surjection onto the set $\{\overline{s}\}_{s\in S}$ which generates $\Sigma(\mathcal{C}(S))$ and so to prove that it is an isomorphism it will suffice to show that it is a monomorphism. Below we discuss injectivity of the map $s\mapsto\overline{s}$ before returning later to a discussion of when the map is a homomorphism, which requires more careful consideration.

The notion of self-automaton semigroups is first introduced by Cain in \cite{cain09}. However, due to the fact that he chooses to act from the right, the class of semigroups obtained by Cain is different from the class introduced here. The precise relationship between these two classes will be discussed later in section \ref{sec:cain}.

First we consider an example of a self-automaton semigroup:
\theoremstyle{definition}
\newtheorem{exln}[def5]{Example}
\begin{exln}
\label{exln}
Let $L_n=\{x_1,x_2,\ldots,x_n\}$ be a left-zero semigroup (i.e. $x_ix_j=x_i$ for all $i,j\in\{1,\ldots,n\}$). Notice that for $i\neq j$ we have $\overline{x_i}\cdot x_i=x_i\neq x_j=\overline{x_j}\cdot x_i$ and hence $\overline{x_i}\neq\overline{x_j}$ showing that the map $x_i\mapsto\overline{x_i}$ is injective. 

For any sequence $\alpha\in L_n^*$ we have $\overline{x_i}\cdot\alpha=(x_i)^k$ where $|\alpha|=k$. Now we have $\overline{x_i}\cdot\overline{x_j}\cdot\alpha=\overline{x_i}\cdot(x_j)^k=(x_i)^k=\overline{x_i}\cdot\alpha=\overline{x_ix_j}\cdot\alpha$ and hence the map $x_i\mapsto\overline{x_i}$ is a homomorphism. Hence we conclude that $L_n$ is self-automaton.
\end{exln}

Motivated by Example \ref{exln} we now move towards establishing when the map $s\mapsto\overline{s}$ is injective in general.
\theoremstyle{definition}
\newtheorem{def6}[def5]{Definition}
\begin{def6}
\label{def6}
Let $S$ be a finite semigroup. For each $a\in S$ define the map $\lambda_a:S\rightarrow S$ by $\lambda_a(x)=ax$ where $x\in S.$ Then $\lambda_a\in\mathcal{T}_S$ and so there is a map $\lambda:S\rightarrow\mathcal{T}_S$ given by $\lambda(a)=\lambda_a.$ The map $\lambda$ is the \emph{left-regular representation} of $S$ and is said to be \emph{faithful} if $\lambda$ is injective.
\end{def6}

The following appears in \cite{maltcev09} but is proved from the perspective of wreath products. We give an alternative proof here by instead acting on sequences:
\theoremstyle{plain}
\newtheorem{cor1}[def5]{Theorem}
\begin{cor1}
\label{cor1}
Let $S$ be a semigroup and let $s,t\in S.$ Then $\overline{s}=\overline{t}\in\Sigma(\mathcal{C}(S))$ if and only if $sa=ta$ for all $a\in S.$
\end{cor1}
\begin{proof}
Let $s,t\in S$ and $\alpha\in S^*.$

($\Rightarrow$) Assume that $\overline{s}=\overline{t}\in\Sigma(\mathcal{C}(S)).$ Let $a\in S.$ Then $\overline{s}\cdot a\alpha=(sa)(\overline{sa}\cdot\alpha)$ and $\overline{t}\cdot a\alpha=(ta)(\overline{ta}\cdot\alpha).$ The outputs must agree on all terms, and in particular the first term. Hence $sa=ta$ for all $a\in S.$

($\Leftarrow$) Assume that $sa=ta$ for all $a\in S.$ Then $\overline{s}\cdot a\alpha=(sa)(\overline{sa}\cdot\alpha)=(ta)(\overline{ta}\cdot\alpha)=\overline{t}\cdot a\alpha$ and hence $\overline{s}=\overline{t}.$
\end{proof}

\newtheorem{th3}[def5]{Lemma}
\begin{th3}
\label{th3}
Let $S$ be a finite semigroup. The map $s\mapsto\overline{s}$ is injective if and only if the left-regular representation of $S$ is faithful.
\end{th3}
\begin{proof}
($\Rightarrow$) Let $x,y\in S$. If $x\neq y$ then $\overline{x}\neq\overline{y}$ and so by Theorem \ref{cor1} there exists $a\in S$ such that $xa\neq ya.$ Hence $$x\neq y\implies\exists a:xa\neq ya\implies\lambda_x\neq\lambda_y\implies\lambda(x)\neq\lambda(y)$$ and the representation is faithful.

($\Leftarrow$) Since $\lambda$ is injective we have, for $x,y\in S$, $$x\neq y\implies\lambda(x)\neq\lambda(y)\implies\lambda_x\neq\lambda_y\implies\exists a:xa\neq ya\implies\overline{x}\neq\overline{y}$$ and hence $s\mapsto\overline{s}$ is injective.
\end{proof}

Having established when $s\mapsto\overline{s}$ is injective, it remains to determine when it is a homomorphism.

At this point, one may pause to consider why we have chosen to define the notion of being self-automaton by using this particular isomorphism, rather than using a direct analogue of the definition given in \cite{cain09} (that is, defining a semigroup to be self-automaton simply when $S\cong\Sigma(\mathcal{C}(S))$ for an arbitrary isomorphism). 

Currently no examples of an isomorphism $\phi:S\rightarrow\Sigma(\mathcal{C}(S))$ have been found where $\phi$ is not the map $s\mapsto\overline{s}$ (we do not consider examples such as $s\mapsto\overline{\theta(s)}$ where $\theta$ is an automorphism of $S$ to be examples of a different isomorphism). Further inspection of the map $s\mapsto\overline{s}$ reveals the following:

\newtheorem{th3a}[def5]{Lemma}
\begin{th3a}
\label{th3a}
Let $S$ be a finite semigroup such that $S\cong\Sigma(\mathcal{C}(S))$. If the map $s\mapsto\overline{s}$ is an injection then it is an isomorphism.
\end{th3a}
\begin{proof}
Let $s,t\in S$. Then there must exist $u\in S$ such that $\overline{s}\cdot\overline{t}=\overline{u}$. By acting on the sequence consisting of a single symbol $a\in S$ we obtain
$$\overline{s}\cdot\overline{t}\cdot a=sta=ua=\overline{u}\cdot a$$
and by Theorem \ref{cor1} we conclude that $\overline{u}=\overline{st}$. Hence the map $s\mapsto\overline{s}$ is an isomorphism.
\end{proof}

An immediate consequence of Lemma \ref{th3a} is that if an example of a semigroup $S$ satisfying $S\cong\Sigma(\mathcal{C}(S))$ but not via $s\mapsto\overline{s}$ exists then it can not have a faithful left-regular representation. Interestingly, there exist examples of semigroups $S$ without faithful left-regular representations satisfying $|S|=|\Sigma(\mathcal{C}(S))|$ but $S\ncong\Sigma(\mathcal{C}(S))$ (such an example is given by the zero-union of a nilpotent monogenic semigroup and a right-zero semigroup of the appropriate size). However, an example where $S\cong\Sigma(\mathcal{C}(S))$ will not be found in this way as if $T$ is a nilpotent semigroup of class $n$ then $\Sigma(\mathcal{C}(T))$ is nilpotent of class $n-1$ (see \cite[Proposition 6.13]{cain09}).

As an example of a semigroup $S$ satisfying $S\cong\Sigma(\mathcal{C}(S))$ but not via $s\mapsto\overline{s}$ appears elusive it would seem that the more restricted Definition \ref{def5} is the correct definition to use (at least for the moment). 

\section{Bands}
\label{sec:bands}

Cain indicates the importance of bands in \cite{cain09} as a source of self-automaton semigroups in his setting. Having already seen Example \ref{exln}, we now show that bands in general provide an abundance of self-automaton semigroups in our setting.

\newtheorem{th4}{Lemma}[section]
\begin{th4}
\label{th4} 
Let $B$ be a finite band. Then the map $b\mapsto\overline{b}$ is a homomorphism.
\end{th4}
\begin{proof}
First notice that for any band $B$ and elements $\beta_1,\beta_2,\ldots,\beta_n\in B$ and for $i\leq j\leq n$ we have $\beta_1\ldots\beta_i\beta_1\ldots\beta_j=\beta_1\ldots\beta_j.$ 

Let $\alpha=\alpha_1\alpha_2\ldots\alpha_n\in B^*.$ Let $s,t\in B.$ We have that
\begin{align*}
\overline{s}\cdot\overline{t}\cdot\alpha&=\overline{s}\cdot(t\alpha_1)(t\alpha_1\alpha_2)\ldots(t\alpha_1\alpha_2\ldots\alpha_n)\\
&=(st\alpha_1)(st\alpha_1t\alpha_1\alpha_2)\ldots(st\alpha_1t\alpha_1\alpha_2\ldots t\alpha_1\ldots\alpha_n)\\
&=(st\alpha_1)(st\alpha_1\alpha_2)\ldots(st\alpha_1\alpha_2\ldots\alpha_n)\\
&=\overline{st}\cdot\alpha.
\end{align*}
Hence $s\mapsto\overline{s}$ is a homomorphism.
\end{proof}

Coupling Lemma \ref{th4} with Theorem \ref{th3} we immediately obtain the following:
\newtheorem{th4a}[th4]{Theorem}
\begin{th4a}
\label{th4a}
A finite band is self-automaton if and only if its left-regular representation is faithful.
\end{th4a}

Having established Theorem \ref{th4a}, one may wonder now what the \textquotedblleft correct analogue\textquotedblright of Cain's open problem should be in terms of left actions. We propose that it should be the following:

\newtheorem{prob}[th4]{Question}
\begin{prob}
\label{prob}
Does the class of self-automaton semigroups consist precisely of those finite bands which have a faithful left-regular representation? 
\end{prob}

In order to answer this, we prove the following:

\newtheorem{th7}[th4]{Theorem}
\begin{th7}
\label{th7}
Let $S$ be a finite semigroup with relative left and right identities (that is, for all $s\in S$ there exist $e,f\in S$ such that $se=fs=s).$ Then $S$ is self-automaton if and only if $S$ is a band with a faithful left-regular representation. 
\end{th7}
\begin{proof}
($\Rightarrow$) Let $s\in S$ and let $e,f\in S$ be such that $se=fs=s.$ Let $\alpha_1,\alpha_2\in S.$ We must have
$$(s\alpha_1)(s\alpha_1 e\alpha_1\alpha_2)=(s\alpha_1)(s\alpha_1 s\alpha_1\alpha_2)$$
since
$$\overline{s}\cdot\overline{e}\cdot\alpha_1\alpha_2=\overline{se}\cdot\alpha_1\alpha_2=\overline{s}\cdot\alpha_1\alpha_2=\overline{fs}\cdot\alpha_1\alpha_2=\overline{f}\cdot\overline{s}\cdot\alpha_1\alpha_2$$
and hence $s\alpha_1 e\alpha_1\alpha_2=s\alpha_1 s\alpha_1\alpha_2$ for all $\alpha_1,\alpha_2\in S.$ By taking $\alpha_1=\alpha_2=e$ we see that $s^2=s$ and hence $S$ is a band.

Since $S$ is self-automaton, the map $s\mapsto\overline{s}$ is injective and so by Lemma \ref{th3} the left-regular representation of $S$ is faithful.

($\Leftarrow$) This follows from Lemmas \ref{th3} and \ref{th4}.
\end{proof}

As a corollary, we deduce that Question \ref{prob} has a positive answer in the following cases:

\newtheorem{th6}[th4]{Corollary}
\begin{th6}
\label{th6}
A finite monoid is self-automaton if and only if it is a band. 
\end{th6}

\newtheorem{th6a}[th4]{Corollary}
\begin{th6a}
\label{th6a}
A finite regular semigroup is self-automaton if and only if it is a band with a faithful left-regular representation.
\end{th6a}

\section{Non-Band Examples}
\label{sec:nonband}

Corollaries \ref{th6} and \ref{th6a} show that Question \ref{prob} has a positive answer in the cases of monoids and regular semigroups. However, we go on to show that the answer in general is negative.

First we prove a result (which is a generalisation of Lemma \ref{th4}) that we will use in Example \ref{ex2}.
\newtheorem{th8}{Lemma}[section]
\begin{th8}
\label{th8}
Let $S$ be a finite semigroup. If $S^2$ is a band then the map $s\mapsto\overline{s}$ is a homomorphism.
\end{th8}
\begin{proof}
First recall that $S^2=\{xy:x,y\in S\}$.

Let $s,t\in S$ and let $\alpha\in S^*.$ Then 
\begin{align*}
\overline{s}\cdot\overline{t}\cdot\alpha&=(st\alpha_1)(st\alpha_1 t\alpha_1\alpha_2)\ldots(st\alpha_1 t\alpha_1\alpha_2\ldots t\alpha_1\ldots\alpha_n)\\
&=(st\alpha_1)(st\alpha_1\alpha_2)\ldots(st\alpha_1\ldots\alpha_n) \text{ since each $t\alpha_1\ldots\alpha_i$ is an idempotent}\\
&=\overline{st}\cdot\alpha.\\
\end{align*}
Hence $s\mapsto\overline{s}$ is a homomorphism.
\end{proof}
\theoremstyle{definition}
\newtheorem{ex2}[th8]{Example}
\begin{ex2}
\label{ex2}
Let $S$ be the semigroup defined by the following Cayley Table:
\begin{center}
\begin{tabular}{c|c c c c}
  & $a$ & $b$ & $c$ & $d$\\
\hline
$a$ & $b$ & $b$ & $b$ & $c$\\
$b$ & $b$ & $b$ & $b$ & $b$\\
$c$ & $c$ & $c$ & $c$ & $c$\\
$d$ & $d$ & $d$ & $d$ & $d$\\
\end{tabular}
\end{center}
Clearly the left-regular representation of $S$ is faithful and so by Lemma \ref{th3} $s\mapsto\overline{s}$ is injective. Observe that $S^2=\{b,c,d\}\cong L_3,$ a three-element left-zero semigroup. Hence by Lemma \ref{th8} $s\mapsto\overline{s}$ is an isomorphism. 
\end{ex2}

This is the first counterexample to Question \ref{prob}.

Next, we exhibit examples of semigroups which satisfy $S=S^2$ and are self-automaton, but which are not bands.
\theoremstyle{definition}
\newtheorem{ex3}[th8]{Example}
\begin{ex3}
\label{ex3}
Let $S_1, S_2,\ldots, S_m$ be finite self-automaton semigroups and define $T=S_1\cup S_2\cup\ldots\cup S_m\cup\{a_{1,1},\ldots,a_{1,n_1},a_{2,1},\ldots a_{2,n_2},\ldots,a_{m,1},\ldots,a_{m,n_m},0\}$ where the product in $S$ extends the products in each $S_i$ and we set $a_{i,j}s_i=a_{i,j}$ for all $j\in\{1,\ldots,n_i\}$, $s_i\in S_i$ and all other products to $0.$

Let $s_{i_1}, s_{i_2}\in S_i.$ Consider the sequence $\alpha=X_1X_2\ldots X_kZB_1B_2\ldots$ where $X_1,\ldots,X_k\in S_i$, $Z\in T\setminus S_i$ and $B_j\in T$. Then
\begin{align*}
\overline{s_{i_1}}\cdot\overline{s_{i_2}}\cdot\alpha&=(\overline{s_{i_1}}\cdot\overline{s_{i_2}}\cdot X_1X_2\ldots X_k)0^\omega\\
&=(\overline{s_{i_1}s_{i_2}}\cdot X_1X_2\ldots X_k)0^\omega \text{ since $S_i$ is self-automaton}\\
&=\overline{s_{i_1}s_{i_2}}\cdot\alpha.\\
\end{align*} 
Notice that by taking the string $X_1X_2\ldots X_k$ to be empty we have accounted for acting on any sequence over $T$ which is not a sequence of elements entirely from $S_i$. If $\beta$ is a sequence of elements entirely from $S_i$ then it follows from the fact that $S_i$ is self-automaton that $\overline{s_{i_1}}\cdot\overline{s_{i_2}}\cdot\beta=\overline{s_{i_1}s_{i_2}}\cdot\beta$. Hence all the products in each $S_i$ hold in $\Sigma(\mathcal{C}(T))$ and so $S_i\leq\Sigma(\mathcal{C}(T))$.

We also have that
\begin{align*}
\overline{a_{i,j}}\cdot\overline{s_{i_1}}\cdot\alpha&=\overline{a_{i,j}}\cdot(\overline{s_{i_1}}\cdot X_1X_2\ldots X_k)0^\omega\\
&=(a_{i,j})^k0^\omega\\
&=\overline{a_{i,j}}\cdot\alpha.\\
\end{align*}
Again, by taking the string $X_1X_2\ldots X_k$ to be empty we have accounted for acting on all sequences over $T$ that are not sequences of elements entirely from $S_i$. If $\beta$ is a sequence of elements entirely from $S_i$ then $\overline{a_{i,j}}\cdot\overline{s_{i_1}}\cdot\beta=(a_{i,j})^{\omega}=\overline{a_{i,j}}\cdot\beta$. Hence $\overline{a_{i,j}}\cdot\overline{s_{i_1}}=\overline{a_{i,j}}$ for all $i\in\{1,\ldots,m\}$ and $j\in\{1,\ldots,n_i\}$ and we conclude that all of these products also hold in $\Sigma(\mathcal{C}(T)).$

Every other product in $T$ is of the form $xy=0$ and so these products will also hold in $\Sigma(\mathcal{C}(T)).$ Hence all products from $T$ hold in $\Sigma(\mathcal{C}(T))$ and so the map $s\mapsto\overline{s}$ is a homomorphism. Using Theorem \ref{cor1}, we show below that the map is also injective:

For $s_1, s_2\in S_i$ there exists $a\in S_i$ such that $s_1a\neq s_2a$ (since $S_i$ is self-automaton) and hence $\overline{s_1}\neq\overline{s_2}.$ For $s_i\in S_i$ and $s_j\in S_j$ ($S_i\neq S_j$) then $s_js_i=0\neq s_is_i\in S_i$ and hence $\overline{s_i}\neq\overline{s_j}.$ If $a_{i,j}\neq a_{k,l}$ then there exists $b\in S$ such that $a_{i,j}b\neq a_{k,l}b$ (we can choose $b\in S_i$). Note also that $s_is_i\neq a_{k,l}s_i$ for all $i,k,l$ and hence $\overline{s_i}\neq\overline{a_{k,l}}.$ Finally observe that for all $x\neq 0$ there exists $y\in S$ such that $xy\neq 0$ and so $\overline{x}\neq\overline{0}$ for all $x\neq 0.$

Hence it is an isomorphism. We have the following egg-box diagram:
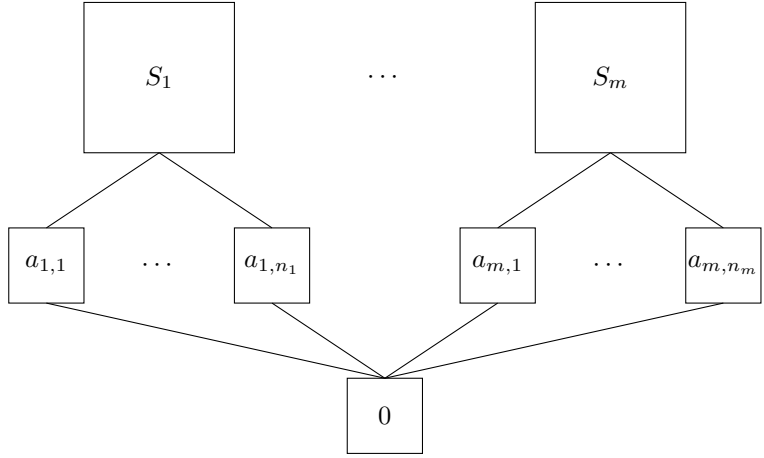
\begin{center}
\begin{tikzpicture}[scale=2]
\draw (-2,2) rectangle (-1,3);
\draw (1,2) rectangle (2,3);
\draw (-2.5,1) rectangle (-2,1.5);
\draw (-1,1) rectangle (-0.5,1.5);
\draw (0.5,1) rectangle (1,1.5);
\draw (2,1) rectangle (2.5,1.5);
\draw (-0.25,0) rectangle (0.25,0.5);
\draw (-2.25,1.5)--(-1.5,2);
\draw (-1.5,2)--(-0.75,1.5);
\draw (0.75,1.5)--(1.5,2);
\draw (1.5,2)--(2.25,1.5);
\draw (-2.25,1)--(0,0.5);
\draw (-0.75,1)--(0,0.5);
\draw (0.75,1)--(0,0.5);
\draw (2.25,1)--(0,0.5);
\draw (-1.5,2.5) node {$S_1$};
\draw (0,2.5) node {$\cdots$};
\draw (1.5,2.5) node {$S_m$};
\draw (-2.25,1.25) node {$a_{1,1}$};
\draw (-0.75,1.25) node {$a_{1,n_1}$};
\draw (0.75,1.25) node {$a_{m,1}$};
\draw (2.25,1.25) node {$a_{m,n_m}$};
\draw (0,0.25) node {$0$};
\draw (-1.5,1.25) node {$\cdots$};
\draw (1.5,1.25) node {$\cdots$};
\end{tikzpicture}
\end{center}
\end{ex3}

\section{Comparisons with Cain's Construction}
\label{sec:cain}

As indicated earlier, we have defined states to act on the left of sequences, in contrast to the approach taken by Cain who acts from the right. The aim of this section is to address the similarities and differences between the two approaches and show how the two are related, before resolving a question stated in the introduction.
\theoremstyle{definition}
\newtheorem{def7}{Definition}[section]
\begin{def7}
\label{def7}
Let $S$ be a finite semigroup. Define $\Pi(\mathcal{C}(S))$ to be the semigroup generated by $\{\overline{s}\}_{s\in S}$ by acting on sequences from the right. That is, for a sequence $\alpha=\alpha_1\alpha_2\ldots\alpha_n$ and states $\overline{s}, \overline{t}$, we have $$\alpha\cdot\overline{s}=(s\alpha_1)(s\alpha_1\alpha_2)\ldots(s\alpha_1\alpha_2\ldots\alpha_n)$$ and $$\alpha\cdot\overline{s}\cdot\overline{t}=(\alpha\cdot\overline{s})\cdot\overline{t}.$$
\end{def7}
\theoremstyle{plain}
\newtheorem{th12}[def7]{Theorem}
\begin{th12}
\label{th12}
Let $S$ be a finite semigroup and $x_1,\ldots,x_n\in S.$ The map $\phi:\Sigma(\mathcal{C}(S))\rightarrow\Pi(\mathcal{C}(S))$ which maps $\overline{x_1}\cdot\ldots\cdot\overline{x_n}\mapsto\overline{x_n}\cdot\ldots\cdot\overline{x_1}$ is an anti-isomorphism.
\end{th12}
\begin{proof}
Recall that a map $\phi$ is an \emph{anti-homomorphism} if, for $x,y\in S$ we have $\phi(xy)=\phi(y)\phi(x).$

In $\Pi(\mathcal{C}(S))$ we have $\phi(\overline{x_1}\cdot\ldots\cdot\overline{x_n})=\overline{x_n}\cdot\ldots\cdot\overline{x_1}=\phi(\overline{x_n})\cdot\ldots\cdot\phi(\overline{x_1}).$ Hence $\phi$ is an anti-homomorphism. The generating sets $\{\overline{s}\}_{s\in S}$ for $\Sigma(\mathcal{C}(S))$ and $\Pi(\mathcal{C}(S))$ are in bijection so we conclude that $\phi$ is an anti-isomorphism.
\end{proof}
By this result, we have that $\Sigma(\mathcal{C}(S))$ and $\Pi(\mathcal{C}(S))$ are dual copies of each other. Hence, to obtain information about $\Pi(\mathcal{C}(S))$ it will suffice to determine $\Sigma(\mathcal{C}(S))$ and then take the dual. Given the way that the action is defined, it is perhaps more natural to work with $\Sigma(\mathcal{C}(S))$ and this is the reason for using left actions throughout this paper, rather than right actions. 

In Cain's setting of right actions, our main question becomes the following: when is the map $S\rightarrow\Pi(\mathcal{C}(S))$ defined by $s\mapsto\overline{s}$ an \emph{anti-isomorphism}? Cain's notion of being self-automaton means that $S\cong\Pi(\mathcal{C}(S))$. Within this section, we refer to such semigroups as being \emph{C-self-automaton}. In the following theorem, we express being C-self-automaton in terms of our setting.
 
\newtheorem{th13}[def7]{Theorem}
\begin{th13}
\label{th13}
Let $s\mapsto\overline{s}$ be an anti-isomorphism $S\rightarrow\Pi(\mathcal{C}(S)).$ Then $S$ is C-self-automaton if and only if $S$ is self-dual and self-automaton.
\end{th13}
\begin{proof}
Recall that a semigroup is said to be \emph{self-dual} if it is anti-isomorphic to itself.

If the map $s\mapsto\overline{s}$ is also an isomorphism then $S$ is commutative and is hence self-dual. By Theorem \ref{th12} $S$ is self-automaton.

If $s\mapsto\overline{s}$ is an anti-isomorphism but not an isomorphism then again by Theorem \ref{th12} $S$ is self-automaton. Suppose that $\phi:\Pi(\mathcal{C}(S))\rightarrow S$ is an isomorphism. Define the map $\psi:S\rightarrow S$ by $\psi(x)=\phi(\overline{x})$ for all $x\in S.$ It is clear that $\psi$ is a bijection.

For $x,y\in S$ we have
$$\psi(xy)=\phi(\overline{xy})=\phi(\overline{y}\cdot\overline{x})=\phi(\overline{y})\phi(\overline{x})=\psi(y)\psi(x).$$
Hence $\psi$ is an anti-isomorphism $S\rightarrow S$ and $S$ is self-dual.

Conversely, it is clear by Theorem \ref{th12} that if $S$ is self-automaton and self-dual then $S\cong\Pi(\mathcal{C}(S)).$
\end{proof}

\newtheorem{rem1}[def7]{Remark}
\begin{rem1}
\label{rem1}
C-self-automaton semigroups are always self-automaton but not conversely. For example, left-zero semigroups (Example \ref{exln}) are self-automaton but are not C-self-automaton as they are not self-dual.
\end{rem1}

In \cite{cain09} Cain asks the question of which semigroups $S$ satisfy $S\cong\Pi(\mathcal{C}(S)).$ Namely,
\newtheorem{conj1}[def7]{Question}
\begin{conj1}
\label{conj1}
Does the class of C-self-automaton semigroups consist of precisely those finite bands in which every $\mathcal{D}$-class is square and every maximal $\mathcal{D}$-class is a singleton?
\end{conj1}
The remainder of the paper will address this question.

\newtheorem{lem1}[def7]{Lemma}
\begin{lem1}
\label{lem1}
Let $S$ be a self-automaton semigroup. If $a,x\in S$ are such that $xa=a$ then $a^2=a.$
\end{lem1}
\begin{proof}
Assume that $S$ is self-automaton. Consider $\overline{x}\cdot\overline{x}$ and $\overline{x^2}$ acting on a sequence $\alpha_1\alpha_2$ and equate the second outputs to obtain $x(x\alpha_1)^2\alpha_2=x^2\alpha_1\alpha_2.$ Setting $\alpha_1=a$ gives $a^2\alpha_2=a\alpha_2$ for all $\alpha_2\in S$ and we conclude by Theorem \ref{cor1} that $a^2=a.$
\end{proof}

We immediately deduce the following as a corollary of Lemma \ref{lem1}.

\newtheorem{lem2}[def7]{Lemma}
\begin{lem2}
\label{lem2}
Let $S$ be a self-automaton semigroup and $a\in S$ be such that $a^2\neq a.$ Then the $\mathcal{L}$-class of $a$ is trivial.
\end{lem2}
\begin{proof}
We know from Lemma \ref{lem1} that if $S$ is self-automaton and $a^2\neq a$ then there does not exist $x\in S$ such that $xa=a.$ Now suppose that $a\mathcal{L}y$ for some $y\in S.$ So there exist $u,v\in S^1$ such that $ua=y$ and $vy=a.$ This gives us $vua=a$ which is a contradiction unless $u=v=1.$ Hence $a=y$ and $L_a$ is trivial.
\end{proof}

\newtheorem{lem3}[def7]{Lemma}
\begin{lem3}
\label{lem3}
Let $S$ be a self-dual, self-automaton semigroup and let $a,x\in S.$ If $ax=a$ then $a$ is an idempotent.
\end{lem3}
\begin{proof}
Let $\phi:S\rightarrow S$ be an anti-isomorphism. Then we have that $\phi(x)\phi(a)=\phi(a)$ and by Lemma \ref{lem1} $\phi(a)$ is an idempotent. Hence $a$ is also an idempotent.
\end{proof}
This means that in self-dual self-automaton semigroups, no non-idempotent elements can be stabilised by multiplication on either side.

\newtheorem{lem4}[def7]{Lemma}
\begin{lem4}
\label{lem4}
Let $S$ be a self-dual, self-automaton semigroup and let $z=xy$ where either $x$ or $y$ is a regular element of $S.$ Then $z^2=z.$
\end{lem4}
\begin{proof}
If $x$ is regular then write $x=qx$ for some $q\in S.$ Then $z=xy=qxy=qz$ and by Lemma \ref{lem1} $z^2=z.$ 

If $y$ is regular then write $y=yp$ for some $p\in S.$ Then $z=xy=xyp=zp$ and by Lemma \ref{lem3} $z^2=z.$
\end{proof}

\newtheorem{th14}[def7]{Theorem}
\begin{th14}
\label{th14}
Let $S$ be a self-automaton and self-dual semigroup. If $S^2=S$ then $S$ is a band. 
\end{th14}
\begin{proof}
Let $a\in S$ and suppose that $a^2\neq a.$ Then by Corollary \ref{th6a} $a$ is not a regular element. We can choose $a$ such that $a$ is in a maximal $\mathcal{D}$-class with respect to the non-regular elements of $S.$ Write $a=bc$ for some $b,c\in S.$ By Lemma \ref{lem4} neither $b$ nor $c$ can be regular elements of $S.$

Since $S$ is self-dual, by Lemma \ref{lem2} $D_a=\{a\}.$ If $b=a$ or $c=a$ then we have either $a=ac$ or $a=ba$ which would be a contradiction by Lemma \ref{lem1} or Lemma \ref{lem3}. Hence $b\neq a$ and $c\neq a.$

This gives us at least 2 non-regular elements in $S.$ Since $a=bc$ we have $D_b>D_a$ which is a contradiction as $D_a$ was assumed to be maximal with respect to the non-regular elements of $S.$ Hence $a^2=a$ and $S$ is a band.
\end{proof}

So we have established that in the case when $S^2=S$ and $s\mapsto\overline{s}$ is an anti-isomorphism it is necessary for $S$ to be a band in order to have $S\cong\Pi(\mathcal{C}(S)).$ Combining this with Theorem \ref{th13} we obtain the following:

\newtheorem{cor2}[def7]{Corollary}
\begin{cor2}
\label{cor2}
The only semigroups satisfying $S^2=S$ and $S\cong\Pi(\mathcal{C}(S))$ (where $s\mapsto\overline{s}$ is an anti-isomorphism)  are the self-dual bands with faithful left-regular representations.
\end{cor2}

If, however, we could find an example of a self-dual semigroup satisfying $S\neq S^2$ and fulfilling the conditions of Lemmas \ref{th3} and \ref{th8}, we would have a counterexample to Question \ref{conj1}. After a discussion of these conditions with Benjamin Steinberg, he suggested the following \cite{steinberg}:
\theoremstyle{definition}
\newtheorem{counterex}[def7]{Example}
\begin{counterex}
\label{counterex}
Let $X=\{1,2,3,4,5\}$ and $X^{\prime}=\{1^{\prime},2^{\prime},3^{\prime},4^{\prime},5^{\prime}\}.$ Let $a,b:X\rightarrow X$ be the functions given by
$$a=\begin{pmatrix}
1 & 2 & 3 & 4 & 5 \\
2 & 3 & 3 & 4 & 5 \\
\end{pmatrix}
b=\begin{pmatrix}
1 & 2 & 3 & 4 & 5 \\   
4 & 5 & 4 & 4 & 5 \\
\end{pmatrix}.
$$
Let $T=\langle a,b\rangle$ where $a,b$ act on the \emph{right} of $X$. We have that $a\neq a^2=a^3$, $b^2=b$ and $ba=b.$ This gives $T=\{a,a^2,b,ab,a^2b\}.$ Note $T\neq T^2=\{a^2,b,ab,a^2b\}$ which is a band.

Now let $a^{\prime}, b^{\prime}:X^{\prime}\rightarrow X^{\prime}$ be given by
$$a^{\prime}=\begin{pmatrix}
1^{\prime} & 2^{\prime} & 3^{\prime} & 4^{\prime} & 5^{\prime} \\
2^{\prime} & 3^{\prime} & 3^{\prime} & 4^{\prime} & 5^{\prime} \\
\end{pmatrix}
b^{\prime}=\begin{pmatrix}
1^{\prime} & 2^{\prime} & 3^{\prime} & 4^{\prime} & 5^{\prime} \\   
4^{\prime} & 5^{\prime} & 4^{\prime} & 4^{\prime} & 5^{\prime} \\
\end{pmatrix}.
$$ 
Let $T^{\prime}=\langle a^{\prime},b^{\prime}\rangle$ where $a^{\prime}, b^{\prime}$ act on the \emph{left} of $X^{\prime}$. We have that $a^{\prime}\neq (a^{\prime})^2=(a^{\prime})^3$, $(b^{\prime})^2=b^{\prime}$ and $a^{\prime}b^{\prime}=b^{\prime}.$ This gives $T^{\prime}=\{a^{\prime},(a^{\prime})^2,b^{\prime},b^{\prime}a^{\prime},b^{\prime}(a^{\prime})^2\}.$ Note $T^{\prime}\neq (T^{\prime})^2=\{(a^{\prime})^2,b^{\prime},b^{\prime}a^{\prime},b^{\prime}(a^{\prime})^2\}$ which is a band.

Note that $T$ and $T^{\prime}$ are dual to each other.

Now let $\hat{T}=\langle(a^{\prime},a),(b^{\prime},b)\rangle\leq T^{\prime}\times T.$ It is easily shown that $|\hat{T}|=11$
and that
\begin{align*}\hat{T}=\{&(a^{\prime},a),(b^{\prime},b),((a^{\prime})^2,a^2),(b^{\prime},ab),(b^{\prime}a^{\prime},b),(b^{\prime},a^2b),(b^{\prime}a^{\prime},ab),\\ &(b^{\prime}(a^{\prime})^2,b),(b^{\prime}a^{\prime},a^2b),(b^{\prime}(a^{\prime})^2,ab),(b^{\prime}(a^{\prime})^2,a^2b)\}.
\end{align*}

Observe that $(\hat{T})^2=\hat{T}\setminus\{(a^{\prime},a)\}$ and $(\hat{T})^2$ is a band.

The egg-box diagram of $\hat{T}$ is:

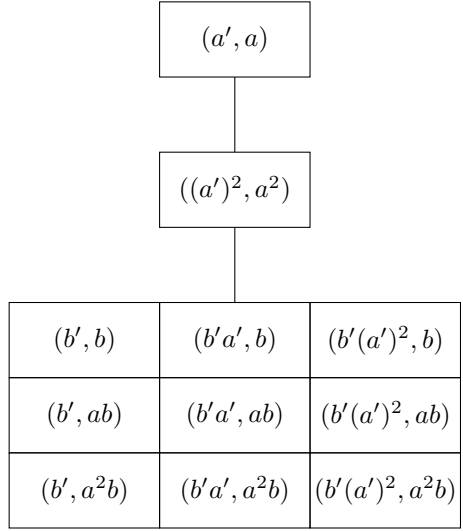
\begin{center}
\begin{tikzpicture}
\draw (-1,6) rectangle (1,7);
\draw (-1,4) rectangle (1,5);
\draw (-3,2) rectangle (-1,3);
\draw (-1,2) rectangle (1,3); 
\draw (1,2) rectangle (3,3);
\draw (-3,1) rectangle (-1,2);
\draw (-1,1) rectangle (1,2);
\draw (1,1) rectangle (3,2);
\draw (-3,0) rectangle (-1,1);
\draw (-1,0) rectangle (1,1);
\draw (1,0) rectangle (3,1);
\draw (0,6)--(0,5);
\draw (0,4)--(0,3);
\draw (0,6.5) node {$(a^{\prime},a)$};
\draw (0,4.5) node {$((a^{\prime})^2,a^2)$};
\draw (-2,2.5) node {$(b^{\prime},b)$};
\draw (0,2.5) node {$(b^{\prime}a^{\prime},b)$};
\draw (2,2.5) node {$(b^{\prime}(a^{\prime})^2,b)$};
\draw (-2,1.5) node {$(b^{\prime},ab)$};
\draw (0,1.5) node {$(b^{\prime}a^{\prime},ab)$};
\draw (2,1.5) node {$(b^{\prime}(a^{\prime})^2,ab)$};
\draw (-2,0.5) node {$(b^{\prime},a^2b)$};
\draw (0,0.5) node {$(b^{\prime}a^{\prime},a^2b)$};
\draw (2,0.5) node {$(b^{\prime}(a^{\prime})^2,a^2b)$};
\end{tikzpicture}
\end{center}

Observe that $\hat{T}$ is self-dual under that map that fixes $(a^{\prime},a), ((a^{\prime})^2,a^2),$ and the diagonal of the bottom $\mathcal{D}$-class, and flips all other elements over the main diagonal. 

The left-regular representation of $\hat{T}$ is not faithful - notice that $\overline{(b^{\prime}a^{\prime},t)}=\overline{(b^{\prime}(a^{\prime})^2,t)}$ for $t\in\{b,ab,a^2b\}.$

To remedy this, let $R=X^{\prime}\times X$, a $5\times 5$ rectangular band. Let $S=\hat{T}\cup R$ with multiplication defined by retaining products from $\hat{T}$ and $R$ and setting
$$(u,v)(i,j)=(u(i),j)$$ and $$(i,j)(u,v)=(i,jv)$$ for $u\in T^{\prime},v\in T,i\in X^{\prime},j\in X.$

It is easily checked that this multiplication is associative.

$R$ is the minimal ideal of $S$ and note that $S^2\neq S$, $S^2$ is a band and $S$ is self-dual. It remains to verify that $S$ has a faithful left-regular representation.

Observe:
$$(b^{\prime}a^{\prime},b)(1,2)=(5,2)\neq(4,2)=(b^{\prime}(a^{\prime})^2)(1,2)$$
and hence $\overline{(b^{\prime}a^{\prime},b)}\neq\overline{(b^{\prime}(a^{\prime})^2,b)}.$

Similarly, $\overline{(b^{\prime}a^{\prime},ab)}\neq\overline{(b^{\prime}(a^{\prime})^2,ab)}$ and $\overline{(b^{\prime}a^{\prime},a^2b)}\neq\overline{(b^{\prime}(a^{\prime})^2,a^2b)}.$

If $i\neq k$ then $(i,j)$ and $(k,l)$ do not act the same on the left of $R$ so $\overline{(i,j)}\neq\overline{(k,l)}.$

If $\{j,k\}\neq\{2,3\}$ then $(i,j)(a^{\prime},a)=(i,ja)\neq(i,ka)=(i,k)(a^{\prime},a)$ and so $\overline{(i,j)}\neq\overline{(i,k)}.$

Also, $(i,2)(b^{\prime},b)=(i,5)\neq(i,4)=(i,3)(b^{\prime},b).$

Hence the left-regular representation of $S$ is faithful. Observe that $|S|=36.$

We have shown that $S$ satisfies the following:
\begin{enumerate}
\item $S\neq S^2$ and hence $S$ is not a band,
\item $S^2$ is a band and so by Lemma \ref{th8} $s\mapsto\overline{s}$ is a homomorphism,
\item $S$ has a faithful left-regular representation,
\item $S$ is self-dual.
\end{enumerate}
Conditions 2 and 3 show us that $S$ is self-automaton. Condition 4 shows that, by Theorem \ref{th12}, $S\cong\Pi(\mathcal{C}(S)).$ Since $S$ is not a band, this is clearly a counterexample to Question \ref{conj1}.
\end{counterex}

\section{Other Results}
\label{sec:other}

We now establish some properties of self-automaton semigroups in general.
\theoremstyle{plain}
\newtheorem{th9}{Lemma}[section]
\begin{th9}
\label{th9}
Let $S$ be a finite semigroup such that $s\mapsto\overline{s}$ is a homomorphism. Let $f:S\rightarrow T$ be an epimorphism of semigroups. Then the map $t\mapsto\overline{t}$ is also a homomorphism $T\rightarrow\Sigma(\mathcal{C}(T))$.
\end{th9}
\begin{proof}
As the map $s\mapsto\overline{s}$ is a homomorphism, notice that, for $x,y\in S$ and $\alpha=\alpha_1\alpha_2\ldots\alpha_n\in S^*$ we have $\overline{x}\cdot\overline{y}\cdot\alpha=\overline{xy}\cdot\alpha$ and hence
$$(xy\alpha_1)\ldots(xy\alpha_1y\alpha_1\alpha_2\ldots y\alpha_1\alpha_2\ldots\alpha_n)=(xy\alpha_1)\ldots(xy\alpha_1\alpha_2\ldots\alpha_n)$$ Hence $xy\alpha_1y\alpha_1\alpha_2\ldots y\alpha_1\alpha_2\ldots\alpha_n=xy\alpha_1\alpha_2\ldots\alpha_n$ for all $n$.

Let $\beta=\beta_1\beta_2\ldots\beta_n\in T^*$ where $\beta_i=f(\alpha_i)$ for some $\alpha_i\in S.$ Then, for $x,y\in S$,
\begin{align*}
\overline{f(x)}\cdot\overline{f(y)}\cdot\beta&=(f(x)f(y)\beta_1)\ldots(f(x)f(y)\beta_1f(y)\beta_1\beta_2\ldots f(y)\beta_1\ldots\beta_n)\\
&=(f(xy\alpha_1))\ldots(f(xy\alpha_1y\alpha_1\alpha_2\ldots y\alpha_1\ldots\alpha_n))\\
&=(f(xy\alpha_1))\ldots (f(xy\alpha_1\alpha_2\ldots\alpha_n))\\
&=(f(xy)\beta_1)\ldots(f(xy)\beta_1\beta_2\ldots\beta_n)\\
&=\overline{f(xy)}\cdot\beta\\
&=\overline{f(x)f(y)}\cdot\beta.\\
\end{align*}
Hence $t\mapsto\overline{t}$ is a homomorphism.
\end{proof}
In particular, if $s\mapsto\overline{s}$ is a homomorphism then it is in fact an epimorphism (since the set $\{\overline{s}\}_{s\in S}$ generates $\Sigma(\mathcal{C}(S))$) and this property is passed to $\Sigma(\mathcal{C}(S)).$ This leads to the following result:
\newtheorem{th10}[th9]{Theorem}
\begin{th10}
\label{th10}
Let $S$ be such that $s\mapsto\overline{s}$ is a homomorphism. Then $\Sigma(\mathcal{C}(S))$ is isomorphic to the image of the  left-regular representation of $S.$
\end{th10}
\begin{proof}
Let $L=\{\lambda_a:x\mapsto ax:a\in S\}$ be the image of the left-regular representation of $S$ and define $\phi:L\rightarrow\Sigma(\mathcal{C}(S))$ by $\phi(\lambda_a)=\overline{a}.$ It is clear that $\phi$ is a surjective map. If $\phi(\lambda_a)=\phi(\lambda_b)$ then $\overline{a}=\overline{b}\implies ax=bx\text{ for all }x\in S\implies\lambda_a=\lambda_b$ and $\phi$ is injective. Hence $\phi$ is a bijection.

We also have $\phi(\lambda_a)\phi(\lambda_b)=\overline{a}\cdot\overline{b}=\overline{ab}=\phi(\lambda_{ab})=\phi(\lambda_a\lambda_b)$ and hence $\phi$ is a homomorphism and $L\cong\Sigma(\mathcal{C}(S)).$ 
\end{proof}

In the case when $S$ is a band, we can go a little further than Theorem \ref{th10}.

\newtheorem{th10a}[th9]{Theorem}
\begin{th10a}
\label{th10a}
Let $S$ be a band. Then $\Sigma(\mathcal{C}(S))$ is self-automaton.
\end{th10a}
\begin{proof}
It suffices to show that $\Sigma(\mathcal{C}(S))$ has a faithful left-regular representation, since by Lemma \ref{th9} we know that the map $\overline{s}\mapsto\overline{(\overline{s})}$ is a homomorphism for all $\overline{s}\in\Sigma(\mathcal{C}(S)).$

Let $a,b\in S$ be such that $\overline{a}\neq\overline{b}.$ Then by Theorem \ref{cor1} there exists $x\in S$ such that $ax\neq bx.$ It now follows by using Lemma \ref{th4} that $axx\neq bxx\implies\overline{ax}\neq\overline{bx}\implies\overline{a}\cdot\overline{x}\neq\overline{b}\cdot\overline{x}$ and hence the left-regular representation is faithful.
\end{proof}

However, it is not possible to replace the hypothesis that $S$ is a band with the assumption that the map $s\mapsto\overline{s}$ is a homomorphism in Theorem \ref{th10a}, as we now show by means of an example:
\theoremstyle{definition}
\newtheorem{ex10b}[th9]{Example}
\begin{ex10b}
\label{ex10b}
Let $S$ be the semigroup defined by the following Cayley Table:
\begin{center}
\begin{tabular}{c|c c c c}
  & $a$ & $b$ & $c$ & $d$\\
\hline
$a$ & $a$ & $b$ & $c$ & $a$\\
$b$ & $a$ & $b$ & $c$ & $a$\\
$c$ & $a$ & $b$ & $c$ & $b$\\
$d$ & $a$ & $b$ & $c$ & $a$\\
\end{tabular}
\end{center}
Then $S^2=\{a,b,c\}\cong R_3$, a three-element right-zero semigroup and so by Lemma \ref{th8} the map $s\mapsto\overline{s}$ is a homomorphism. Note also that $\overline{a}=\overline{b}=\overline{d}$ by Theorem \ref{cor1}. 

It now follows that $\overline{x}\cdot\overline{y}=\overline{xy}=\overline{y}$ where $x,y\in\{a,c\}$ and hence $\Sigma(\mathcal{C}(S))\cong R_2$ which does not have a faithful left-regular representation so is not self-automaton. In fact, $\Sigma(\mathcal{C}(R_2))\cong\{1\}.$
\end{ex10b}
It is also worth noting that self-automaton semigroups are closed under taking direct products.
\theoremstyle{plain}
\newtheorem{th11}[th9]{Theorem}
\begin{th11}
\label{th11}
$S,T$ are self-automaton semigroups if and only if $S\times T$ is also self-automaton.
\end{th11}
\begin{proof}
In order to prove this, we use a result due to Mintz \cite{mintz09} which states that for finite semigroups $S$ and $T$, $\Sigma(\mathcal{C}(S\times T))\leq\Sigma(\mathcal{C}(S))\times\Sigma(\mathcal{C}(T))$ and $\overline{(s,t)}$ can be written as $(\overline{s},\overline{t})$ as $\Sigma(\mathcal{C}(S\times T))$ is isomorphic to the subsemigroup of $\Sigma(\mathcal{C}(S))\times\Sigma(\mathcal{C}(T))$ where the words in each component have the same length.

($\Rightarrow$) Let $S,T$ be self-automaton semigroups. Define $\phi:S\times T\rightarrow \Sigma(\mathcal{C}(S\times T))$ by $\phi((s,t))=\overline{(s,t)}.$ This map is clearly a surjection to the generating set $\{\overline{(s,t)}\}_{s\in S, t\in T}$. If $\overline{(s,t)}=\overline{(u,v)}$ then by Theorem \ref{cor1} $(s,t)(\alpha,\beta)=(u,v)(\alpha,\beta)$ for all $\alpha\in S, \beta\in T\implies (s\alpha,t\beta)=(u\alpha,v\beta)\implies s\alpha=u\alpha, t\beta=v\beta\implies s=u$ and $t=v$ since $S,T$ are self-automaton $\implies(s,t)=(u,v)$ and the map is injective.

We also have that 
\begin{align*}
\overline{(s,t)}\cdot\overline{(u,v)}&=(\overline{s},\overline{t})\cdot(\overline{u},\overline{v})\\
&=(\overline{s}\cdot\overline{u},\overline{t}\cdot\overline{v})\\
&=(\overline{su},\overline{tv}) \text{ since $S$ and $T$ are self-automaton}\\
&=\overline{(su,tv)}\\
&=\overline{(s,t)(u,v)}\\
\end{align*}
and hence $\phi$ is an isomorphism.

($\Leftarrow$) Assume that $S\times T$ is self-automaton. Then $\overline{(s,t)}\cdot\overline{(u,v)}=\overline{(s,t)(u,v)}=\overline{(su,tv)}.$ Hence by equating the outputs of $\overline{(s,t)}\cdot\overline{(u,v)}\cdot\gamma$ and $\overline{(su,tv)}\cdot\gamma$ (where $\gamma=(\alpha_1,\beta_1)(\alpha_2,\beta_2)\ldots$ and $\alpha_i\in S, \beta_j\in T$) we obtain 
$$su\alpha_1u\alpha_1\alpha_2\ldots u\alpha_1\ldots\alpha_n=su\alpha_1\ldots\alpha_n$$
and
$$tv\beta_1v\beta_1\beta_2\ldots v\beta_1\ldots\beta_n=tv\beta_1\ldots\beta_n$$
for all $n.$ This forces $\overline{s}\cdot\overline{u}\cdot\alpha=\overline{su}\cdot\alpha$ and $\overline{t}\cdot\overline{v}\cdot\beta=\overline{tv}\cdot\beta$ (where $\alpha=\alpha_1\alpha_2\ldots$ and $\beta=\beta_1\beta_2\ldots$). Hence the maps $s\mapsto\overline{s}$ and $t\mapsto\overline{t}$ are homomorphisms.

If $s_1\neq s_2$ then $\overline{(s_1,t)}\neq\overline{(s_2,t)}$ for all $t\in T.$ Hence by Theorem \ref{cor1} there exists $(a,b)\in S\times T$ such that $(s_1,t)(a,b)\neq(s_2,t)(a,b)\implies(s_1a,tb)\neq(s_2a,tb)\implies s_1a\neq s_2a\implies\overline{s_1}\neq\overline{s_2}$ and hence $s\mapsto\overline{s}$ is injective. Similarly, $t\mapsto\overline{t}$ is injective.

It now follows that $s\mapsto\overline{s}$ and $t\mapsto\overline{t}$ are isomorphisms.
\end{proof}

\section*{Acknowledgements}
The author is extremely grateful to Nik Ru\v{s}kuc and Collin Bleak for their many helpful contributions and valuable discussions whilst this work was being undertaken.

The final publication is available at Springer via\\
http://dx.doi.org/10.1007/s00233-014-9610-3.

\end{document}